\documentclass[12pt]{article}
\usepackage{graphicx}%
\usepackage{multirow}%
\usepackage{amsmath,amssymb,amsfonts}%
\newcommand{\norm}[1]{\left\lVert#1\right\rVert}
\usepackage{amsthm}%
\newtheorem{prop}{Proposition}
\newtheorem{lem}{Lemma}
\newtheorem{thm}{Theorem}
\usepackage{mathrsfs}%
\usepackage[title]{appendix}%
\usepackage{xcolor}%
\usepackage{textcomp}%
\usepackage{manyfoot}%
\usepackage{booktabs}%
\usepackage{algorithm}%
\usepackage{algorithmicx}%
\usepackage{algpseudocode}%
\usepackage{listings}%
\usepackage{verbatim}
\usepackage{physics}
\usepackage{hyphenat}
\usepackage[square,numbers]{natbib}
\usepackage{authblk}
\newenvironment{keywords}{%
   \small	
   \quotation
   \textbf{Keywords:}
   \endquotation
}

\begin{document}

\title{Asymptotic expansions for approximate solutions of boundary integral equations
}
\author{Akshay Rane, Kunalkumar Shelar} 
\affil{Institute of Chemical Technology, Mumbai, India}

\date{January 2024}
\maketitle

\begin{abstract}
This paper uses the Modified Projection Method to examine the errors in solving the boundary integral equation from Laplace's equation. The analysis uses weighted norms, and parallel algorithms help solve the independent linear systems. By applying the method developed by Kulkarni, the study shows how the approximate solution behaves in polygonal domains. It also explores computational techniques using the double-layer potential kernel to solve Laplace’s equation in these domains. The iterated Galerkin method provides an approximation of order $2r+2$ in smooth domains. However, the corners in polygonal domains cause singularities that reduce the accuracy. By adjusting the mesh near these corners, the accuracy can almost be restored when the error is measured using the uniform norm. This paper builds on the work of Rude et al. By using Kulkarni's modified operator and observes superconvergence in iterated solutions. This leads to an asymptotic error expansion, with the leading term being $O(h^4)$ and the remaining error term $O(h^6)$, resulting in a method with similar accuracy.
\end{abstract}
AMS Subject classification 45B05,65R20
\footnote{First author would like to thank to UGC-FRP for their support.}


\section{Introduction}
The boundary integral equation (BIE) method has emerged as a powerful approach for solving potential problems, especially those governed by Laplace's equation. One of its main advantages is that it reduces the dimensionality of the problem by focusing on the boundary, resulting in fewer unknowns and often higher computational efficiency compared to traditional finite element or finite difference methods. This benefit is particularly pronounced when dealing with complex geometries, such as polygonal domains, which introduce additional challenges like boundary singularities at corners. While the BIE method offers improved efficiency with fewer unknowns, handling these singularities is crucial for maintaining convergence rates and accuracy.

Existing approaches, such as the Galerkin method with piecewise polynomials, provide effective solutions for smooth domains. However, in polygonal domains, the presence of corners degrades the convergence rate due to singularities in the solution and kernel. Mesh grading techniques have been introduced by Rice \cite{RICE1976} to mitigate this issue, restoring much of the convergence when measured in the uniform norm. For comprehensive information on the numerical treatment of boundary integral equations, one can refer to \cite{Atkinson1990BOUNDARYIE} and Hackbusch \cite{Hackbusch1991}.
Chandler \cite{Chandler_1984} demonstrated that the iterated Galerkin solution can achieve superconvergence, with accuracy improvements up to order 2 in the uniform norm. Furthermore, Richardson extrapolation has been shown to increase the convergence order to 4, as detailed by Rude et al., who also introduced multi-parameter extrapolation for handling boundary integral equations in polygonal domains. \\
Lin and Xie \cite{Lin_Xie} demonstrated the existence of asymptotic error expansions under suitable conditions. This paper builds on the work of Rude et al.\cite{Rude} by extending their multi-parameter asymptotic expansion technique using a modified operator proposed by Kulkarni \cite{Kulkarni2009ExtrapolationUA}. This finite rank approximating operator enhances the accuracy of solving second-kind Fredholm integral equations, allowing for superconvergence effects in iterated solutions. The resulting asymptotic error expansion achieves a leading term of \(O(h^4)\), where \(h\) is the mesh characteristic parameter, with the remaining error term being \(O(h^6)\). By using linear combinations of solutions for different mesh parameters, we can eliminate the leading error term, achieving a highly accurate solution with a significantly reduced computational cost. 
In addition, we explore multi-parameter asymptotic expansions where the boundary is partitioned into multiple segments, each with an independently chosen mesh width. This technique results in a discretization scheme with several independent parameters, offering enhanced flexibility in mesh refinement and accuracy. Our work extends previous analyses, providing new insights and computational strategies for improving the accuracy and efficiency of boundary integral equation solutions, particularly for complex, polygonal geometries.
In multi-parameter asymptotic expansions, the boundary is partitioned into  $r$ parts. The mesh width $h_i$ for each segment $(1 \leq i \leq r)$ can be chosen independently of the others.
 This results
in a discretization that has $r$ independent parameters $(h_1,..., h_r)$. We will extend the analysis of \cite{Lin_Xie} to show that the resulting discrete approximation permits a multivariate asymptotic error expansion in these p parameters.

In section 2, we state the precise
definition of our problem set notations, and recall some preliminaries, In section 3 we describe the modified projection method.
Asymptotic series expansion for the proposed solution is obtained in
Section 4.
Section 5 summarizes our results. It also contains implementation details and a comparison with the method suggested by Rude et al \cite{Rude}.

\section{Preliminaries}

\noindent Let $\Omega$ be a simply connected polygon in $\mathbb{R}^2 $ with boundary $\bar{\Gamma}$ and corner points $x_0,x_1,...,x_r$ where $x_r=x_0$.

\noindent The corresponding interior angle at $x_i$ is denoted by $(1-p_i)\pi$, for $i=0,1,...r$  for $| {p_i}|<1$ and $ p_r=p_0$. 
Consider double layered potential formulation of Laplace's equation in $\Omega$ leading to a boundary integral equation of the second kind on $\bar{\Gamma}$:

\begin{equation}\label{eq:1}
u(x)+Tu(x)=f(x), x\in \bar{\Gamma}
\end{equation}

where
\begin{equation}\label{eq:2}
\begin{split}
    Tu(x)=\int_{\bar{\Gamma}} k(x,y)u(y)ds_y, \\k(x,y)=\frac{1}{\pi} \frac{\partial}{\partial n_y} ln|x-y|
\end{split}
\end{equation}

\begin{figure}[ht]
    \centering
    \includegraphics[width=0.8\linewidth]{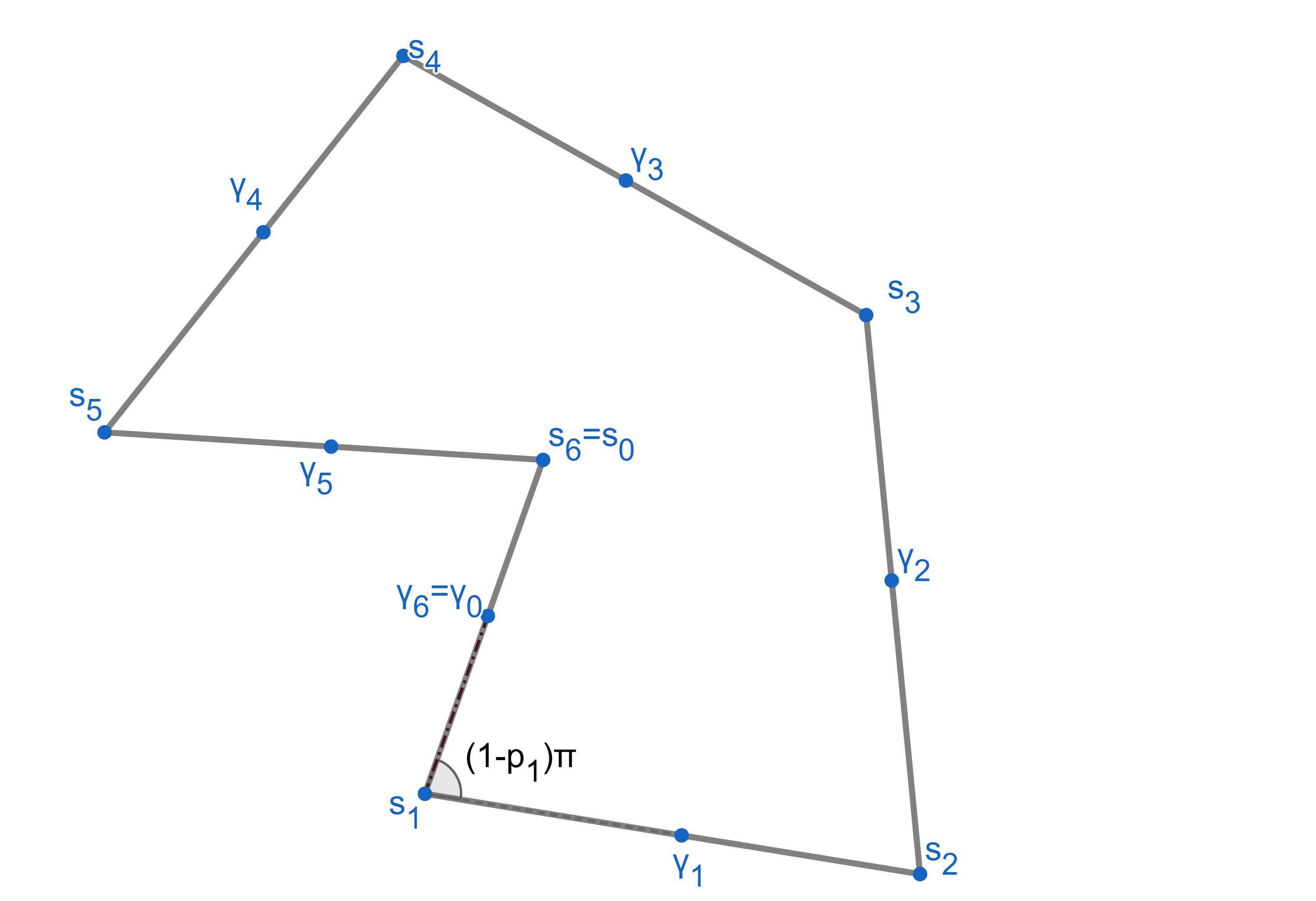}
    \caption{\small{Example of Polygonal Domain and initial partitioning of the boundary}}
    \label{fig:enter-label}
\end{figure}

\noindent Let $ \bar{\Gamma}=x(s)$  be the parametrization, where the boundary be parameterized by the arc length s for $s_0 \leq s<s_r $, with $s=s_j$ as the corner points corresponding to $x_j$ for $j=0,...,r.$ We will not distinguish the arc length parameter s from the corresponding point $x=x(s)$ on $\bar{\Gamma}$ when the meaning is clear from the context.
Assume that $\Gamma_p=\{{\bar{\Gamma}_1,...,\bar{\Gamma}_r}\}$ is an initial partition of $\bar{\Gamma}$ such that the corner point $x_j$ is an inner point of $\bar{\Gamma}_j.$
Here $\bar{\Gamma}_j=x(\Gamma_j),$ with 
$\Gamma_j=[\gamma_{j-1},\gamma_j]$ and
 $\gamma_{j-1}<s_j<\gamma_j$, for $j=0,...,r-1.$ \\For $j=r$, we define $\Gamma_r=[\gamma_{r-1},s_r] \cup [s_0,\gamma_0].$\\
 See figure \ref{fig:enter-label}.

\noindent Values $s_j$ and $\gamma_j$ follows strictly increasing sequence
$s_0<\gamma_0<s_1<\gamma_1<s_2<...<\gamma_r-1<s_r$. We introduce $\gamma_r=\gamma_0.$
\\
The kernel on $\Gamma_j$ is the well-known expression

\begin{equation*}
 \eval{k(s,t)}_{\Gamma_j \times \Gamma_j} =
\left\{
	\begin{array}{ll}
		 \frac{\sin (p_j)}{\pi} \frac{|s-s_j|}{(s-s_j)^2+(s_j-t)^2+2(s-s_j)(s_j-t)\cos (p_j \pi)} & \mbox{if }  s<s_j<t \text{ or } t<s_j<s  \\
		0 & \mbox{if } s_j<s,t \text{ or } s,t<s_j<s
	\end{array}
\right.
\end{equation*}

\noindent We introduce weighted norms. For $\alpha_j>0$ and an integer $m\geq 0$, we first define for each part $\Gamma_j$ of the boundary 

\begin{equation*}
\norm{u}_{m,\alpha_j}=\max _{0\leq i \leq m} \sup\{|s-s_j|^{max\{ i-\alpha_j,0\}}|D^i u(s)|:s \in \Gamma_j\backslash \{s_j\}\}
\end{equation*}

\begin{equation*}
    C_{\alpha_j}^m(\Gamma_j)=\{u\in C^m(\Gamma_j \backslash \{s_j\}) : \norm{u}_{m,\alpha_j}<\infty\}
\end{equation*}

\noindent and use this norm to construct a norm for function on $\Gamma$ by

\begin{equation*}
 \norm{u}_{m,\alpha}=\max_{1\leq j\leq r} \norm{w\restriction_{\Gamma_j}}_{m,\alpha_j}
\end{equation*}

\begin{equation*}
 C_{\alpha}^m(\Gamma)=\bigcap_{j=1}^{p} C_{\alpha_j}^m(\Gamma_j)
\end{equation*}
\noindent where $\alpha=(\alpha_1,...,\alpha_r)$ is a vector of parameters.
$(C_{\alpha}^m(\Gamma), \norm{.}_{m,\alpha})$ is a Banach space.
Our construction has partitioned the domain $\Gamma$ into 2r parts, given by $\Gamma_j^+=[s_j,\gamma_j]$ and $\Gamma_j^-=[\gamma_{j-1}, s_j]$. Thus $u=(u^1,...,u^{2r})$ on $\Gamma$ can be written as
\begin{equation*}
    u^{2j-1}=\eval{u}_{\Gamma_j^-} \hspace{12pt} and 
    \hspace{12pt}
    u^{2j}=\eval{u}_{\Gamma_j^+}
\end{equation*}

\noindent Consider now
\begin{equation*}
   T_j ^{+}u_j^+=\int_{\Gamma_j^+} k(s,t)u_j^+ dt ,\hfill s\in \Gamma_j^-, 
\end{equation*}
    
\begin{equation*}
   T_j ^{-}u_j^-=\int_{\Gamma_j^-} k(s,t)u_j^- dt , s\in \Gamma_j^+ .
\end{equation*}
on $\Gamma_j$, the operator can now be represented in $2\times 2$ block as 
\begin{equation}\label{eq:3}
    T_j=\begin{pmatrix}
    0 & T_j^{+}\\
    T_j ^{-}&0
\end{pmatrix}
\end{equation}

\noindent with the above notation , the integral equation \eqref{eq:1} is equivalent to a $2r\times 2r$ system of equations
\begin{equation}\label{eq:4}
    (I+\tilde{T})\tilde{u}=\tilde{f}
\end{equation}
\\
with $\tilde{u}=(u^1,...u^{2r})$ and $\tilde{f}=(f^1,...f^{2r})$.Here $\tilde{T}$ is a matrix that can be decomposed into
\begin{equation}\label{eq:5}
    \tilde{T}=diag[T_1,...T_r]+S,
\end{equation}
where S is a system with smooth component kernels.

\noindent For the sake of simplicity, we will refer  $\tilde{T}$  in  \eqref{eq:2} as T, since they are equivalent matrix representations. Similarly, $\tilde{u}$ and $\tilde{f}$ will be denoted as $u$ and $f$, respectively, henceforth.

\noindent We write the equation as
\begin{equation}\label{main_eq}
    (I+T)u=f
\end{equation}

\noindent From identity
\begin{equation*}
    \int_{0}^{\infty}\frac{y^\beta}{x^2+y^2+2xy\cos\alpha}=\frac{\pi \sin(\beta \alpha)}{\sin (\alpha) \sin(\beta \pi)}, |\beta|<1, |\alpha|<\pi ,
\end{equation*}
The operator $R:L^p(0,\infty)\rightarrow L^p(0,\infty) $ where $ 1<p\leq\infty$ defined by
\begin{equation*}
    Ru(x)=\frac{\sin{(p_j \pi)}}{\pi}\int_{0}^{\infty} \frac{x u(y)dy}{x^2+y^2+2xy} \cos{(p_j\pi)}
\end{equation*}
is linear, bounded and
\begin{equation*}
    \norm{R}_{L^p\rightarrow L^p}\leq\frac{\sin(|v_j|\pi /p)}{\sin (\pi/ p)}.
\end{equation*}
Thus, as an operator from $L^p(\Gamma_j)$ into $L^p(\Gamma_j)$, we have
\begin{equation*}
    \norm{T_j}_{L^p\rightarrow L^p}\leq\frac{\sin(|p_j|\pi/ p)}{\sin (\pi /p)}
\end{equation*}

\noindent For  $ p\rightarrow \infty $, we get
\begin{equation*}
    \norm{T_j}_{L^\infty\rightarrow L^\infty}\leq |p_j| <1,
\end{equation*}
which gives the bounded system.
\begin{equation*}
    \norm{(I+diag[T_1,...T_r])^{-1}}_{L^\infty\rightarrow L^\infty}<\infty
\end{equation*}
The inverse of the Fredholm operator $I+T$ is also bounded.

\begin{equation}   \label{boundedop}
\norm{(I+T)^{-1}}_{L^\infty\rightarrow L^\infty}<\infty,
\end{equation}

\noindent Define 
$\alpha^{*} \equiv (\alpha_1^*,...\alpha_r^*)$ as 
$\alpha_j^*=\frac{1}{(1+|p_j|)}, 1\leq j\leq r.$

\noindent For all $\alpha < \alpha^{*}$,
$f \in {C_{\alpha}^{m}} $ implies $ u \in {C_{\alpha}^{m} } $ and hence 

\begin{equation}
    \norm{(I+T)^{-1}}_{C_{\alpha}^{m}\rightarrow C_{\alpha}^{m}}<\infty,
\end{equation}

\section{Modified Projection Method}
Consider the equation (\ref{main_eq})
\begin{equation*}
    (I+T)u=f    
\end{equation*}
$(I+T)$ is invertible and bounded as shown in (\ref{boundedop}).
T is compact and $\pi_h \rightarrow I $ pointwise. 
Consider the finite rank operator $T_n^M$ to approximate T. This operator proposed in \cite{Kulkarni2009ExtrapolationUA}.
$$T_n^M=\pi_hT\pi_h+\pi_hT(I-\pi_h)+(I-\pi_h)T\pi_h $$
 The operator $T_n^M$ is a finite linear combination of bounded operators. The sum of bounded operators is also a bounded operator. Therefore, $T_n^M$ is well-defined as a bounded operator on the function space.
 We can see that $n\rightarrow \infty$

$$\norm{T-T_n^M}=\norm{(I-\pi_h)T(I-\pi_h)}\rightarrow 0$$
If we assume the given equation 
has a unique solution then for all large enough n and $u_n^M\in S^h$ 
\begin{equation}
    (I+T_n^M)u_n^M=f    
\end{equation}
also has a unique solution.

\noindent An iterated approximation is given by
\begin{equation}
    \Tilde{u}_n^M=f-Tu_n^M
\end{equation}
It is shown in \cite{kulkarni2003superconvergence} that under specific conditions $u_n^M$ converges more rapidly to $u$ compared to approximations obtained in the Galerkin and Sloan methods. Iterated approximation $\Tilde{u}_n^M$ converges faster to $u$ than $u_n^M$.
We will show that 
\begin{equation*}
    \norm{u-u_n^M}_{\infty}=\mathcal{O}(h^3)
\end{equation*}
and the error in the iterated solution is
\begin{equation*}
    \norm{u-\tilde{u}_n^M}_{\infty}=\mathcal{O}(h^4)
\end{equation*}

\noindent Note that $\Tilde{u}_n^M=\Tilde{u}(h_1,h_2,...h_r)$ is obtained from the discretization with parameters $(h_1,h_2,...h_r)$  for the corresponding segments of $\Gamma$, which we will see in the next section.

\noindent We get
\begin{equation}
    \begin{split}
        u-u_n^M& =((I+T)^{-1}-(I+T_n^M)^{-1})f \\
        &=(I+T_n^M)^{-1}((I+T_n^M)(I-T)^{-1}-I)f\\
        &=(I+T_n^M)^{-1}(T_n^M-T)u
    \end{split}
\end{equation}

\begin{equation}
\begin{split}
    \Tilde{u}_n^M-u=T(u-u_n^M)
    &= T(I+T)^{-1}(T-T_n^M)(I+T_n^M)^{-1}f\\
    &=(I+T)^{-1}T(I-\pi_h)T(I-\pi_h)(u+u_n^M-u)
\end{split}   
\end{equation}
We know $\pi_h$ is uniformly bounded and From \cite{Rude} $$\norm{T(I-\pi_h)}= \mathcal{O}(h^2)$$
and we will show that $$\norm{T(I-\pi_h)T(I-\pi_h)}= \mathcal{O}(h^4)$$

\section{Discretization}
We will now discuss the discretization of equation \eqref{eq:1} for numerical solution.

\noindent Corresponding to each element $\Gamma_j$ of the boundary $\Gamma$ define the mesh size as $h_j$, calculated by $h_j = n_j^{-1}$. Here, $n_j$ represents the number of intervals in the $j^{th}$ segment.The parameter $q_j$ indicates the grading strength towards the $j^{th}$ corner.\\
For each element $\Gamma_j$ of the boundary $\Gamma$, the mesh consists of $2n_j$ intervals.

\begin{equation*}
        \Gamma_j^{h_j}=\{\tau=[t_{j,i},t_{j,i+1}]:0\leq i\leq 2n_j -1\},
    \end{equation*}
 where the nodes $t_{j,i}$ for $i=0,1,...2n_j $ are given by
 \begin{equation*}
     t_{j,i}=s_j - \left(\frac{n_j -i}{n_j}\right)^{q_{j}} (s_j-\gamma_{j-1}), 0\leq i\leq n_j,
 \end{equation*}
 \begin{equation*}
     t_{j,n_j+i}=s_j + \left(\frac{i}{n_j}\right)^{q_{j}} (\gamma_{j}-s_j), 0\leq i\leq n_j,
 \end{equation*}
Here we can observe the points are denser near the corner in every element. The grading exponents $q_j$'s are real numbers with $q_j\geq 1$. 

\noindent For $j=1,...,r,$ the intervals $\Gamma_j^{h_j}$ form a partition
\begin{equation*}
    \Gamma^h=\bigcup_{j=1}^{r}\Gamma_j^{h_j}
\end{equation*}

 \noindent of $\Gamma$, where  
 \begin{equation*}
 h=\max_{1\leq j\leq r}{h_j}.    
 \end{equation*}

\noindent We introduce the space of piece-wise constant functions on $\Gamma^h$ to develop a Galerkin discretization.
 
 $$S^h=\{ v\in L^\infty(\Gamma):  {v} \text { is a constant for all } \tau \in \Gamma^h \}$$

\noindent $S^h$ is piecewise constant polynomials with breakpoints ${t_{j,i}}$. Hence discontinuities occur at breakpoints. As basis functions of $S^h$ are discontinuous,  thus $\pi_h$ and $(I-\pi_h)$ are bounded.\\
Consider $h_j=t_{j,i}-t_{j,i-1}$ for the interval $\Gamma_j^{h_j}=[t_{j,i},t_{j,i+1}]$.\\
Let $\pi_h$ denote the orthogonal projection of $L^2(\Gamma)$ onto space of piecewise constant functions $S^h$.
\\ for any interval $\tau$ in $T^h$ we have 
 $\pi_h\equiv \pi_{h_1,...,h_r}:L^2(\Gamma)\rightarrow S^h$ is defined by
 \begin{equation*}
     \pi_h u=\frac{1}{|\tau|}\int_{\tau}u ds 
 \end{equation*}
\\
We can check $\norm{\pi_h}_{\infty} \leq 1$. Chandler \cite{Chandler_1984} [Lemma 1,2] has shown
$$
\max   \biggl\{ \norm{T\pi_h}_\infty, \norm{\pi_h T}_\infty \biggr \}\leq \abs{p_i}<1    
$$

\noindent Which gives $(I+\pi_hT)^{-1}$ and $(I+T\pi_h)^{-1}$ are bounded.

\noindent We have
 \begin{equation*}
     \norm{(I+\pi_hT)^{-1}}_{L^\infty\rightarrow L^\infty} <\infty
 \end{equation*}
 
  \begin{equation*}
     \norm{(I+T\pi_h)^{-1}}_{L^\infty\rightarrow L^\infty} <\infty
 \end{equation*}

\noindent Assume $u_h \in S^h $ be the Galerkin approximation of (\ref{main_eq}).
$$ (I+\pi_hT)u_h=\pi_hf $$
\noindent The \textit{iterated Galerkin approximation} is given by 
$$u_h^*=f-Tu_h$$ 
where $u_h=u(h_1, ..., h_r)$ 
represents a function $u$ discretized with parameters $h_1,h_2,...,h_r$  for corresponding segments of a domain $\Gamma$ i.e. $h_j \in \Gamma_j^{h_j}$. We can write $$ u_h^*-u=T(u-u_h)  \text{ and } (I+T\pi_h)=T(I-\pi_h)u$$
\noindent Hence error of iterated approximation can be given as
\begin{equation}
    u_h^*-u=(I+T\pi_h)^{-1} T(I-\pi_h)u
\end{equation}
Consider $u_0=u(h_1,h_2,.., h_j,..,h_r)$ be a Galerkin solution over a given mesh. We refine the given mesh by dividing $j^{th}$ element by 2 to get a solution $u_j=u(h_1,h_2,..,h_{j-1}, h_j/2,h_{j+1},..,h_r)$. \\
Therefore $u_0^*$ and $u_j^*$ represent corresponding iterated solutions.\\
\noindent By Chandler's theorem \cite{Chandler} we have \\
$$ 
\norm{u_{h}^{*}-u}_{\infty}\leq ch^2 
$$
\noindent We recall the following proposition \cite{Rude} \\
\begin{prop}
    If $ 0<\alpha<\alpha^{*}, u \in C_\alpha^1{(\Gamma)}, \beta=(\beta_1,...,\beta_r) $ satisfying $\beta_j=\alpha_j-2/q_j$ and $q_j>2/\alpha_j$ then 
\begin{equation}
\norm {T(I-\pi_h)u}_{0,\infty}\leq Ch^2 \norm{u}_{1,\alpha}
\end{equation}
and
\begin{equation}
\norm {T(I-\pi_h)u}_{1,\beta}\leq Ch^2 \norm{u}_{1,\alpha}
\end{equation}
\end{prop} 
\noindent Recalling the following lemma \cite{Rude}
\begin{lem}{\label{lem 3.2}}
    If $u\in C^4(\tau), g\in C^3(\tau),$ then
    \begin{equation*}
        \int_\tau ((u-\pi_h u)g)(t) dt=\frac{1}{12}h_\tau^2\int_\tau Du(t)Dg(t)dt+\gamma_\tau,\\
        \text{ where } h_\tau=|\tau|
        \text{ and }
            \end{equation*}
        \begin{equation*}
        |\gamma_\tau|\leq Ch^4_\tau\sum_{i=0}^{3} \norm{D^{4-i}u}_{0,\infty,\tau} \norm{D^i g}_{0,1,\tau}.
    \end{equation*}    
\end{lem}

\noindent Rude et. al \cite{Rude} had shown following asymptotic expansion exists.

\begin{prop}
\label{prop 3.3}
If $ u \in C_\alpha^4{(\Gamma)}$ , $ q_j >4/\alpha_j, j=1,2,...r$ then 
\begin{equation}
(T(I-\pi_h)u)(s)=\frac{1}{12} \sum_{j=1}^{r} \sum_{\tau\in T^h, s_j \notin \tau} h_\tau^2 \int \partial_t {k(s,t)} Du(t) dt + O(h^4) \norm{u}_{4,\alpha_j} 
\end{equation}
\end{prop}

\noindent We will use proposition 2 to show the following result.
\begin{prop}  \label{prop}
If $u\in C^4_\alpha(\Gamma)$, $q_j>4/\alpha_j$ for $j=1,2,...r$ 
$\alpha_j=a_j-2/q_j$
and $a_j>2/q_j$
then   
\begin{equation}
    \begin{split}
        &(T(I-\pi_h)T(I-\pi_h)u)(s)\\
        &=\frac{1}{144} \sum_{j_1=1}^r \sum_{j_2=1}^r \sum_{\tau_1\in T_1^h, s_j\notin \tau_1} \sum_{\tau_2\in T_2^h, s_j\notin \tau_1} 
     h_{\tau^*}^{4}  \int_{\tau_1}\partial_t k(s,t) \left( \int_{\tau_2}\partial_t\partial_s k(s,t) Du(t)dt\right)dt
      + \mathcal{O}(h^6_j)\norm{u}_{8,\alpha_j} 
    \end{split}
\end{equation}

\noindent Proof:
\end{prop}

\noindent We have 
\begin{equation}
    \begin{split}
    ( T_j^+(I-\pi_h)u)(s)&=\int_{\Gamma_{j}^+}k(s,t)((I-\pi_h)u)(t)dt\\
    &=\frac{1}{12} \sum_{i=n_j+1}^{2n_j-1}h_{\tau_{j,i}}^2 \int_{\tau_{j,i}} \partial_t k(s,t) Du(t)dt + \mathcal{O}(h^4_j)\norm{u}_{4,\alpha_j}
\end{split}
\end{equation}

\noindent Consider
\begin{equation*}
    \phi(u(t))=(T_j^+(I-\pi_h)u)(t)
\end{equation*}
 
\noindent We get
\begin{equation}\label{phi}
 \begin{split}
     (T_j^+(I-\pi_h)\phi(u))(s) 
      & = {\frac{1}{12}} \sum_{i=n_j+1}^{2n_j-1} h_{\tau_{j,i}}^2 \int_{\tau_{j,i}} \partial_t k(s,t)D \phi (u(t))dt + 
    \mathcal{O}(h^4_j)\norm{\phi(u)}_{4,\alpha_j}
\end{split}   
\end{equation}

\begin{equation*}
    D\phi(u(t))=D( T_j^+(I-\pi_h)u)(s)=    \int_{\Gamma_{j}^{+}} \partial_{s}k(s,t)((I-P_h)u)(t)dt
\end{equation*}
\noindent and,

\begin{equation}\label{derphi}
     \sum_{i=n_j+1}^{2n_j-1} \int_{\tau_{j,i}}\partial_{s}k(s,t)((I-\pi_h)u)(t) dt 
    = \frac{1}{12}  \sum_{i=n_j+1}^{2n_j-1} h_{\tau_{j,i}}^{2} \int_{\tau_{j,i}}\partial_{t}\partial_{s}k(s,t)Du(t) dt + \gamma_{\tau_{j,i}} 
\end{equation}
\noindent from lemma (\ref{lem 3.2})
    $$ \lvert \gamma _{\tau_{j,i}} \rvert \leq Ch_{j,i}^4 \sum_{l=0}^3\norm{D^{4-l}u}_{0,\infty,\tau{j,i}}\norm{\partial_t\partial_s k(.,.)}_{0,1,\tau_{j,i}}\\
$$
\noindent Note that
\begin{equation*}
    \begin{split}
        \lvert \gamma _{\tau_{j,i}} \rvert &\leq Ch_{j}^4 (t_{j,i}-s_j)^{4-4/q_j} \norm{u}_{4,\alpha_j} \\ &\hspace{5mm}\cross \sum_{l=0}^3 (t_{j,i}-s_j)^{\alpha_j+l-4} \int_{\tau_{j,i}}((s_j-s)^2 + (t-s_j)^2)^{-(l+3)/2}dt \\
        &\leq C.h_{j}^4 (t_{j,i}-s_j)^{\alpha_j-4/q_j}    \int_{\tau_{j,i}} ((s_j-s)^2 + (t-s_j)^2)^{-3/2}dt \norm{u}_{4,\alpha_j}\\
        &\leq Ch_{j}^4 \norm{u}_{4,\alpha_j} \int_{\tau_{j,i}} ((s_j-s)^2 + (t-s_j)^2)^{-3/2+ \alpha_j/2 -2/q_j}dt
    \end{split}
\end{equation*}
\noindent we obtain 
\begin{equation}\label{tau}
    \sum_{i=n_j+1}^{2n_j-1} \lvert \gamma _{\tau_{j,i}} \rvert \leq Ch_{j}^4 \norm{u}_{4,\alpha_j}
\end{equation}
\noindent provided $q_j > 4/\alpha_j$
\begin{equation*}
    \sum_{i=n_j+1}^{2n_j-1} \int_{\tau_{j,i}}\partial_{s}k(s,t)((I-\pi_h)u)(t)dt 
    = \frac{1}{12}  \sum_{i=n_j+1}^{2n_j-1} h_{\tau_{j,i}}^{2} \int_{\tau_{j,i}}\partial_{t}\partial_{s}k(s,t)Du(t)dt + \mathcal{O}(h^4_j)\norm{u}_{4,\alpha_j}  
\end{equation*}

\noindent Note that
\begin{equation*}
    \int_{\tau_{j,n_j}}\partial_{s}k(s,t)((I-\pi_h)u)(t)dt \leq Ch^4_j\norm{u}_{4,\alpha_j}
\end{equation*}
\noindent So we get bound for (\ref{derphi})\\

\noindent Consider
\begin{equation*}
    \begin{split}
        \norm{\phi(u)}_{4,\alpha} &=\underset{0\leq i\leq4}{max} sup \{|s-s_j|^{max(i-\alpha_j,0)} |D^i\phi(u)(s)|\colon s \in \Gamma_j-{s_j}\}
    \end{split}
\end{equation*}

\noindent We will use the orthogonality property of the projection.
\begin{equation*}
\begin{split}
    \int_{\tau_{j,i}} \partial_s k(s,t)((I-\pi_h)u)(t)dt &=\int_{\tau_{j,i}} (I-\pi_h)\partial_s k(s,t)((I-\pi_h)u)(t)dt \\ 
\end{split}   
\end{equation*}

\begin{equation*}
    \begin{split}
       &=\int_{\tau_{j,i}} (I-\pi_h)(I-\pi_h)\partial_s k(s,t)((I-\pi_h)(I-\pi_h)u)(t)dt\\
    &\leq Ch_{\tau_{j,i}}^2 \norm{\partial_t\partial_t\partial_s k(s,t)}_{0,1,\tau_{j,i}}\norm{D^2u}_{0,\infty,\tau_{j,i}}\\
    &\leq Ch_j^2(t_{j,i}-s_j)^{2-2/q_j} (t_{j,i}-s_j)^{\alpha_j-2}\int_{\tau_{j,i}}\biggl((s-s_j)^2+(t-s_j)^2 \biggr)^{-5/2}dt \norm{u}_{2,\alpha_j}\\
    &\leq Ch_j^2\int_{\tau_{j,i}}\biggl((s-s_j)^2+(t-s_j)^2 \biggr)^{\alpha_j/2-1/q_j-5/2}dt \norm{u}_{2,\alpha_j}  
    \end{split}   
\end{equation*}

\noindent Thus the sum is bounded by
\begin{equation*}
    \begin{split}
      \lvert \sum_{i=n_j+1}^{2n_j-1} \int_{\tau_{j,i}} \partial_s k(s,t)((I-\pi_h)u)(t)dt \rvert &\leq Ch_j^2\int_{\Gamma_j^+}\biggl((s-s_j)^2+(t-s_j)^2 \biggr)^{\alpha_j/2-1/q_j-5/2}dt \norm{u}_{2,\alpha_j}\\
         &\leq Ch_j^2(s-s_j)^{\alpha_j-2/q_j-4}dt \norm{u}_{2,\alpha_j}\\
        &=Ch_j^2(s-s_j)^{\beta_j-4}dt \norm{u}_{2,\alpha_j}
    \end{split}
\end{equation*}
\noindent Where $\beta_j=\alpha_j-2/q_j$
and $\alpha_j>2/q_j$
Also $u\in C^2_\alpha(\Gamma)$.

\noindent We can generalize the above result to get the following:
\begin{equation*}
    \begin{split}
        \int_{\tau_{j,i}} \partial_s k(s,t)((I-\pi_h)u)(t)dt &=\int_{\tau_{j,i}} (I-\pi_h)\partial_s k(s,t)((I-\pi_h)u)(t)dt        
\end{split}   
\end{equation*}

\begin{equation*}
    \begin{split}
       &=\int_{\tau_{j,i}} (I-\pi_h)(I-\pi_h)\partial_s k(s,t)((I-\pi_h)(I-\pi_h)u)(t)dt \\
       & =\int_{\tau_{j,i}} (I-\pi_h)(I-\pi_h)(I-\pi_h)(I-\pi_h)\partial_s k(s,t) ((I-\pi_h)(I-\pi_h)(I-\pi_h)(I-\pi_h)u)(t)dt 
  \end{split}
\end{equation*}
\begin{equation*}
    \begin{split}
    &\leq Ch_{\tau_{j,i}}^2 \norm{\partial_t\partial_t\partial_t\partial_t\partial_s k(s,t)}_{0,1,\tau_{j,i}}\norm{D^4u}_{0,\infty,\tau_{j,i}}\\
    &\leq Ch_j^2(t_{j,i}-s_j)^{2-2/q_j} (t_{j,i}-s_j)^{\alpha_j-4}\int_{\tau_{j,i}}\biggl((s-s_j)^2+(t-s_j)^2 \biggr)^{-9/2}dt \norm{u}_{4,\alpha_j}\\
    &\leq Ch_j^2\int_{\tau_{j,i}}\biggl((s-s_j)^2+(t-s_j)^2 \biggr)^{\alpha_j/2-1/q_j-11/2}dt \norm{u}_{4,\alpha_j} 
    \end{split}
\end{equation*}

\noindent Thus the sum is bounded by
\begin{equation*}
    \begin{split}
         &\abs {\sum_{i=n_j+1}^{2n_j-1}
        \int_{\tau_{j,i}} \partial_s k(s,t)((I-\pi_h)u)(t)dt } \\
       &\leq Ch_j^2\int_{\Gamma_j^+}\biggl((s-s_j)^2+(t-s_j)^2 \biggr)^{\alpha_j/2-1/q_j-11/2}dt \norm{u}_{4,\alpha_j}\\
        &\leq Ch_j^2(s-s_j)^{\alpha_j-2/q_j-10}dt \norm{u}_{4,\alpha_j}\\
        &=Ch_j^2(s-s_j)^{\beta_j-10}dt \norm{u}_{4,\alpha_j}
    \end{split}
\end{equation*}
\noindent Where $\beta_j=\alpha_j-2/q_j$
and $\alpha_j>2/q_j$
Also $u\in C^4_\alpha(\Gamma)$.
\noindent Further, we show the result follows for $i=n_j$.
Hence we get
\begin{equation}
    \norm{T_j^+(I-\pi_h)u}_{4,\beta_j,\tau_{i,j}}\leq Ch_j^2 \norm{u}_{4,\alpha_j}
\end{equation}
Using analogous statements we can prove 
\begin{equation}
    \norm{T_j^+(I-\pi_h)u}_{0,\infty,\Gamma_j^-}\leq Ch_j^2 \norm{u}_{4,\alpha_j}
\end{equation}
Note for
\begin{equation*}
    \norm{\phi(u)}_{4,\alpha_j,\Gamma_j^-} =\norm{K_j^+(I-\pi_h)u}_{4,\alpha_j,\Gamma_j^-}\leq Ch_j^2 \norm{u}_{4,a_j}
\end{equation*}
\noindent where $\alpha_j=a_j-2/q_j$ and $a_j>2/q_j$ 
\begin{equation}\label{normphi}
    \norm{\phi(u)}_{1,\beta_j,\Gamma_j^-} =\norm{T_j^+(I-\pi_h)u}_{1,\beta_j,\Gamma_j^-}\leq Ch_j^2 \norm{u}_{4,\alpha_j}\leq Ch_j^2
\end{equation}

\noindent Hence, equations
\eqref{phi},\eqref{derphi},\eqref{normphi}  give us
$$
            (T_j^+(I-\pi_h)T_j^+(I-\pi_h)u)(s)=T_j^+(I-\pi_h)\phi(u)  
$$
\begin{multline*}
     = {\frac{1}{12}} \sum_{i=n_{j_1}+1}^{2n_{j_1}-1}h_{\tau_{j_1,i_1}}^2 \int_{\tau_{j_1,i_1}} \partial_t k(s,t)\left ( \frac{1}{12}  \sum_{i=n_{j_2}+1}^{2n_{j_2}-1} h_{\tau_{j_2,i_2}}^{2} \int_{\tau_{j_2,i_2}}\partial_{t}\partial_{s}k(s,t)Du(t)dt +\gamma_{\tau}\right )dt \\
        +\mathcal{O} (h^4_j)Ch_j^2\norm{u}_{4,a_j}
\end{multline*}    
\begin{equation*}
    \begin{split}
            =\frac{1}{144} \sum_{i={{n_j}_1}+1}^{{{2n_j}_1}-1}\sum_{i={{n_j}_2}+1}^{{{2n_j}_2}-1}{h_\tau}_{j_1 i_1}^2{h_\tau}_{j_2 i_2}^2\int_{\tau_{j_1 i_1}} \partial_t k(s,t)\left( \int_{\tau_{j_2,i_2}}\partial_{t}\partial_{s}k(s,t)Dudt\right ) dt \\+ {\frac{1}{12}} \sum_{i=n_{j_1}+1}^{2n_{j_1}-1}h_{\tau_{j_1,i_1}}^2 \int_{\tau_{j_1,i_1}} \partial_t k(s,t) \gamma_{\tau} dt+  \mathcal{O}(h^6_j) \norm{u}_{4,a_j}
    \end{split}
\end{equation*}

\begin{equation}
    \begin{split}
      &=\frac{1}{144} \sum_{i={n_j}_1+1}^{{{2n_j}_1}-1}\sum_{i={n_j}_2+1}^{{{2n_j}_2}-1}{h_\tau}_{j_1 i_1}^2{h_\tau}_{j_2 i_2}^2\int_{\tau_{j_1 i_1}} \partial_t k(s,t)
      \left( \int_{\tau_{j_2 i_2}}\partial_{t}\partial_{s}k(s,t)Dudt\right ) dt+ \phi_1+ \mathcal{O}(h^6_j) \norm{u}_{4,a_j}\\
      &=\mathcal{O}(h^4)+ \mathcal{O}(h^6)+ \mathcal{O}(h^6_j) \norm{u}_{4,a_j}\\
        &=\mathcal{O}(h^4)+ \mathcal{O}(h^6) \norm{u}_{4,a_j}  
\end{split}
\end{equation}
\noindent where from (\ref{tau})
$$
       |\phi_1|=\frac{1}{12}\sum_{i={{n_j}_1}+1}^{{2n_j}_1-1}  
       {{h^2}_\tau}_{j_1 i_1} 
       \int_{\tau_{j_1 i_1}}  |\partial_t k(s,t)|\gamma_{\tau} dt\leq Ch^6
$$ 
\noindent We have a relation 
$
    \norm{(I-\pi_h)u}_{0,\infty,\tau_{j,n_j}}\leq Ch^{\alpha_j}_{\tau_{j,n_j}} \norm{u}_{1,\alpha_j}
$ This gives us 
$
    \norm{(T-T_n^M)u}_{0,\infty,\tau_{j,n_j} }=\norm{(I-\pi_h)T(I-\pi_h)u}_{0,\infty,\tau_{j,n_j}} \leq Ch^{\alpha_j}_{\tau_{j,n_j}}\norm{T(I-\pi_h)u}_{1,\alpha_j} 
    =\mathcal{O}(h^3)
$
\noindent we get
\begin{equation}
    \begin{split}
       \norm{u-u_n^M}_{\infty}& =\norm{(I+T_n^M)^{-1}(T_n^M-T)u}_{\infty}\\
       & \leq\norm{(I+T_n^m)^{-1}}_{\infty}\norm{(T_n^m-T)u}_{\infty}\\
      & \leq C \norm{(I-\pi_h)T(I-\pi_h)u}_{\infty}\\
      &=\mathcal{O}(h^3)
    \end{split}
\end{equation}

\noindent Further, we get
\begin{equation}\label{iter1}
    \begin{split}
    \Tilde{u}_n^M-u &=(I+T)^{-1}T(I-\pi_h)T(I-\pi_h)(u+u_n^M-u)\\
    &= {(I+T)^{-1}} ( T(I-\pi_h)T(I-\pi_h)u + (I+T)^{-1} T(I-\pi_h)T(I-\pi_h) ({u_n^M-u}) \\
    &= {(I+T)^{-1}}w_h + {(I+T)^{-1}} 
    T(I-\pi_h)T(I-\pi_h) (I+T_n^M)^{-1}(T_n^M-T)u\\
    &= {(I+T)^{-1}}w_h + {(I+T)^{-1}} 
    T(I-\pi_h)T(I-\pi_h) (I+T_n^M)^{-1}(I-\pi_h)T(I-\pi_h)u\\
    &={(I+T)^{-1}}w_h + {(I+T)^{-1}}(I+T_n^M)^{-1} 
    T(I-\pi_h)T(I-\pi_h)T(I-\pi_h)u\\
    &={(I+T)^{-1}}w_h + {(I+T)^{-1}}(I+T_n^M)^{-1} 
    T(I-\pi_h)w_h\\
    \end{split}
\end{equation}
\noindent where $w_h\equiv w(h_1,h_2,...,h_r)=T(I-\pi_h)T(I-\pi_h)u$
\begin{equation} \label{bound1}
    \begin{split}
        \norm{{(I+T)^{-1}} (I+T_n^M)^{-1}T(I-\pi_h)w_h}_{0,\infty}
        & \leq C \norm{T(I-\pi_h)w_h}\\
        & \leq c h^2 \norm{w_h}_{1,\beta}\\
        & \leq c h^6 \norm{u}_{1,a}
    \end{split}
\end{equation}

\noindent Note:\\
    Let $(I+T)^{-1}=L \Rightarrow$
      $T=(I+T)L=L+TL$\\
      $L=T(I-L) \Rightarrow L\in Range(T)$

\noindent This shows $(I+T)^{-1}T  \in T$.  \\
\noindent  $(I+T)^{-1}T$ is an integral operator and the kernel of this operator will have the same properties as that of kernel of T.Hence without loss of generaality in estimation $(I+T)^{-1}T$ can be treated as T. 
      
\noindent Now (\ref{iter1}) and (\ref{bound1}) lead to

\begin{equation}\label{bound extrapol}
    \Tilde{u}_n^M-u ={(I+T)^{-1}}w_h+ \mathcal{O}(h^6) \norm{u}_{1,a}
\end{equation}
thus proposition (\ref{prop}) produces
\begin{multline}\label{multextra}
        \norm{\frac{1}{3}\biggl(4 \sum_{j=1}^r u(h_1,..., h_{j-1},h_j/2,h_{j+1},...,h_r)-(4r-3)u(h_1,...,h_r)}_{0,\infty} =\mathcal{O}(h^6) \norm{u}_{1,a}
\end{multline}
combining (\ref{bound extrapol}) and (\ref{multextra}) gives us following theorem.

\begin{thm}\label{thmmultetra}
    \begin{multline}
    If u \in C^{4}_{\alpha}(\Gamma), \alpha < \alpha^{*}, \alpha_j=2-q_j, q_j>2/a_j, for j=1,2,...,r ,then \\ \norm{\frac{1}{3}\biggl(4 \sum_{j=1}^r u_j^*-(4r-3)u_0^*}_{0,\infty} \leq Ch^6   
    \end{multline}
 \end{thm}

\section{Conclusion}
since we have proved existance of asymptomatics expansion of approximate solution.Richardson extrapolation can be used to improve the order of convergence of approximate solution.
    
\bibliography{ref}
    
\end{document}